# The Joint Distance Measure: A Measure of Similarity Accounting for Spatial and Angular Distances

Abeeb A. Awotunde

Dept. of Petroleum Engineering and CIPR

King Fahd University of Petroleum & Minerals, Dhahran, Saudi Arabia

**Abstract**

Vector similarity measures play a fundamental role in various fields, including machine learning, natural language processing, information retrieval, and data mining. These measures quantify the closeness between two vectors in a high-dimensional space and are vital for tasks such as document similarity, recommendation systems, and clustering. While several vector similarity measures exist, each similarity measure is suited to specific purposes. In addition, some of these measures lack robustness to certain data types and/or probability distributions. This article presents a measure of vector similarity, known as the Joint Distance Measure (JDM). The JDM combines the Minkowski distance measure and the cosine similarity measure to establish a distance measure that captures the advantages of both measures.

## 1. Introduction

In many applications, data can be represented as vectors in a high-dimensional space, where each vector represents an entity such as a document, user behavior, or image. Measuring the similarity or distance between these vectors is crucial for various algorithms that rely on proximity-based methods. Vector similarity measures help determine the degree of similarity between vectors, allowing systems to group similar entities or rank them by relevance. For instance, similarity measures find widespread application in recommender systems. Recommender systems [1, 2, 3, 4, 5] are algorithms used to predict and suggest items to users based on their preferences, behavior, or interactions with the system. They are commonly used in online platforms such as e-commerce sites (like Amazon), streaming services (like Netflix), social media platforms (like Facebook), and even news websites. These systems aim to enhance user experience by providing personalized recommendations, making it easier for users to discover new content or products they might like. Measures such as cosine similarity and Pearson correlation coefficient [6] help in finding similar users or items to make recommendations.

Another area in which similarity measures have become very useful is natural language processing [7, 8, 9, 10, 11]. Natural Language Processing (NLP) is a field of artificial intelligence (AI) focused on enabling computers to understand, interpret, and generate human language in a way that is both meaningful and useful. The goal of NLP is to allow machines to interact with humans using natural language, such as speech or text, in a way that mimics human understanding and reasoning. In NLP, word embeddings like Word2Vec [12, 13, 14], GloVe [15], and BERT [16, 17, 18] represent words as vectors. Vector similarity measures are then used to calculate word or sentence similarity, which is essential in tasks such as semantic similarity, clustering, and translation.

A third area in which similarity measures have found relevance is information retrieval (IR). Information Retrieval [19, 20, 21, 22, 23] is the process of obtaining relevant information from a large collection of unstructured data, typically stored in databases or search engines, based on a user's query or information needs. It plays a critical role in systems like search engines (Google, Bing), digital libraries, and recommendation systems. The goal of IR is to return the most relevant documents or data items that satisfy a user's request. In information retrieval, vector similarity measures are used to assess the relevance of documents with respect to a given query. Cosine similarity is commonly used to rank documents based on how similar they are to the query vector.

Similarity measures are also used in clustering and classification. Clustering [24, 25, 26, 27] involves grouping similar data points together into clusters, where the data points within a cluster are more similar to one another than to those in other clusters. The goal of clustering is to identify patterns in data without using any labels or predefined categories. Clustering algorithms like K-means clustering [28] and hierarchical clustering [29] rely on distance or similarity measures to group similar data points. Classification [30, 31, 32, 33] is a supervised learning technique. It involves training a model on labeled data, where each data point is associated with a predefined label. The goal of classification is to predict the label or category of a new, unseen data point based on its features. Similarity measures are also used in classification tasks to identify the most relevant features for decision-making.

In this article, we propose the Joint Distance Measure (JDM), a similarity measure that combines the Minkowski distance measure (MDM) with the cosine similarity measure (CSM). The JDM is a more general measure than each of MDM and CSM in that it provides the advantages of both measures, accounting for spatial distance (in the Minkowski's sense) and angular difference. Furthermore, the JDM has better capabilities to distinguish between pairs of vectors than does either MDM or CSM.

In the remainder of this paper, we provide a description of vector similarity measures that are relevant to this article, their mathematical foundations, practical applications, and challenges. In particular, we discuss the cosine similarity, dot product measure, Euclidean distance and city block distance in relation to the development of the joint distance measure. We also present examples that help us study the behavior, validity and accuracy of these similarity measures.

## 2. Mathematical Background

A vector $\vec{u}$ in an $n$-dimensional space can be represented as:

$$\vec{u} = [u_1, u_2, \ldots, u_n]. \tag{1}$$

For two vectors $\vec{u}$ and $\vec{v}$, a similarity measure $S(\vec{u},\vec{v})$ or a distance measure $D(\vec{u},\vec{v})$ quantifies how close or similar these vectors are. Given the vectors $\vec{u}$, $\vec{v}$ and $\vec{w}$, a similarity or distance measure is said to be a metric if it satisfies the following axioms:

1. The distance between a vector and itself is zero: $D(\vec{u},\vec{u})=0$.
2. The distance between any two distinct vectors is positive: If $\vec{u} \neq \vec{v}$, then $D(\vec{u},\vec{v})>0$.
3. Symmetry is established such that the distance from $\vec{u}$ to $\vec{v}$ is the same as the distance from $\vec{v}$ to $\vec{u}$: $D(\vec{u},\vec{v})=D(\vec{v},\vec{u})$.
4. The triangle inequality holds: $D(\vec{u},\vec{w}) \leq D(\vec{u},\vec{v})+D(\vec{v},\vec{w})$.

The last axiom (Axiom 4) states that the distance from a vector $\vec{u}$ to another vector $\vec{w}$ cannot be larger than the sum of the distance between $\vec{u}$ and $\vec{v}$, and that between $\vec{v}$ and $\vec{w}$.

The following are commonly used similarity/distance measures:

---

## 2.1 Cosine Similarity Measure (CSM)

Cosine similarity measure [34, 35, 36, 37, 38] is one of the most widely used measures, especially in text mining and natural language processing. It measures the cosine of the angle between two vectors in a vector space. Thus, CSM is a direct indication of the angular difference between any two vectors. The CSM is given by:

$$S_{\cos}(\vec{u},\vec{v}) = \frac{\vec{u} \cdot \vec{v}}{\|\vec{u}\|\|\vec{v}\|}, \tag{2}$$

where $\vec{u} \cdot \vec{v}$ is the dot product between the vectors $\vec{u}$ and $\vec{v}$, $\|\vec{u}\|$ and $\|\vec{v}\|$ are the magnitudes of $\vec{u}$ and $\vec{v}$, respectively. Alternatively, we may define CSM as the cosine of the angle between the two vectors so that:

$$S_{\cos}(\vec{u},\vec{v}) = \cos\theta, \tag{3}$$

where $\theta$ is the angle between vectors $\vec{u}$ and $\vec{v}$. CSM ranges from $-1$ (completely dissimilar) to 1 (completely similar), with 0 indicating no similarity (decorrelation). A pair of parallel vectors pointing in the same direction have a CSM value of 1, two parallel vectors pointing in opposite directions have a CSM value of $-1$ while a pair of orthogonal vectors have a CSM of 0. Thus, under this measure, a higher similarity index indicates a higher level of similarity between the vectors. CSM reflects the angular difference between two vectors but gives no indication of the spatial distance between the vectors. This is the main deficiency of the cosine similarity measure. Consequently, CSM yields the same value in a case where two vectors make the same angle with an index vector regardless of the spatial distances of these vectors to the index vector. For example, two parallel vectors pointing in the same direction will have the same CSM with respect to a third vector (the index vector) regardless of whether the parallel vectors are identical or not. For certain tasks, particularly where the Euclidean distance between vectors are important, this deficiency can be a major concern. Nonetheless, the CSM has found very useful applications in tasks such as NLP, information retrieval, recommender systems and text mining.

**2.2 Dot Product Similarity Measure (DPSM)**

The DPSM [39, 40, 41, 42] between any two vectors $\vec{u}$ and $\vec{v}$ is the dot product between the vectors, given by:

$$S_{dp}(\vec{u}, \vec{v}) = \vec{u} \cdot \vec{v}. \tag{4}$$

In summation notation, Eq. 4 may be written as:

$$S_{dp}(\vec{u}, \vec{v}) = \sum_{i=1}^{n} u_i v_i, \tag{5}$$

where $u_i$ and $v_i$ are the components at index $i$ of vectors $\vec{u}$ and $\vec{v}$, respectively.

The magnitude and sign of dot product *similarity* measure (DPSM) is said to indicate how similar the vectors are, with increasing values indicating high levels of similarity [39, 40, 41]. In [42], it was mentioned that the DPSM is used in some fields, including natural language processing, and recommendation systems, to quantify the similarity between two vectors. However, we found this measure to be an inconsistent measure of similarity and we do not recommend its use in any task. To explain why the DPSM is inconsistent, first we write the dot product between the vectors $\vec{u}$ and $\vec{v}$ in another form as

$$S_{dp}(\vec{u}, \vec{v}) = \|\vec{u}\| \cdot \|\vec{v}\| \cos\theta. \tag{6}$$

This shows that the dot product is a product of the magnitudes of the two vectors being compared and the cosine of the angle between the vectors. Now consider that we would like to find the similarities of two vectors $\vec{u}_1$ and $\vec{u}_2$ to an index vector $\vec{u}_0$. Then,

$$S_{dp}(\vec{u}_0, \vec{u}_1) = \|\vec{u}_0\| \cdot \|\vec{u}_1\| \cos\theta_{0,1}, \tag{7}$$

and

$$S_{dp}(\vec{u}_0, \vec{u}_2) = \|\vec{u}_0\| \cdot \|\vec{u}_2\| \cos\theta_{0,2}, \tag{8}$$

where $\theta_{0,1}$ is the angle between $\vec{u}_0$ and $\vec{u}_1$, and $\theta_{0,2}$ is the angle between $\vec{u}_0$ and $\vec{u}_2$. We observe that a common factor in Eqs. 7 and 8 is $\|\vec{u}_0\|$. Because this term has equal contributions to the right-hand-side (RHS) of both equations, it may be factored out. Hence, the comparison is between the terms $\|\vec{u}_1\|\cos\theta_{0,1}$ and $\|\vec{u}_2\|\cos\theta_{0,2}$ in Eqs. 7 and 8, respectively. From these two terms, we see that the only portions that make any comparison to the index vector are $\cos\theta_{0,1}$ and $\cos\theta_{0,2}$, which are the respective cosine similarity indices of the two vectors to the index vector. $\|\vec{u}_1\|$ and $\|\vec{u}_2\|$ are magnitudes that scale these similarity indices. This scaling of the similarity indices may distort the true nature of the similarity thereby giving a false notion of similarity. For instance, if the angles $\theta_{0,1}$ and $\theta_{0,2}$ are both acute angles and $\|\vec{u}_2\|$ is much larger than $\|\vec{u}_1\|$, then $S_{dp}(\vec{u}_0,\vec{u}_2)$ will be greater than $S_{dp}(\vec{u}_0,\vec{u}_1)$ and the DPSM will automatically indicate that $\vec{u}_2$ is more similar to $\vec{u}_0$ than $\vec{u}_1$ is. A more vivid illustration of the inconsistency of this measure of similarity is presented in the example section.

## 2.3 Euclidean Distance Measure (EDM)

Euclidean distance [37, 38, 43, 44] is the straight-line distance between two vectors $\vec{u}$ and $\vec{v}$ in Euclidean space. It is the $l_2$-norm of the difference between two vectors and it is given by:

$$D_{Euc}(\vec{u},\vec{v}) = \|\vec{u}-\vec{v}\|_2, \tag{9}$$

where $\|\vec{u}-\vec{v}\|_2$ represents the $l_2$-norm of $\vec{u}-\vec{v}$. Equation 9 may also be written as:

$$D_{Euc}(\vec{u},\vec{v}) = \sum_{i=1}^{n}(u_i-v_i)^2, \tag{10}$$

where $u_i$ and $v_i$ are the values on the $i^{th}$ dimension of the vectors $\vec{u}$ and $\vec{v}$, respectively. The EDM is always positive and a lower EDM value indicates a higher level of similarity between the pair of vectors being compared. The EDM reflects a direct comparison of the individual dimensions of the vectors being compared. Thus, if the differences between corresponding elements in the two vectors are small, the Euclidean distance will be small. In addition, the Euclidean distance is sensitive to the magnitude of the difference of the vectors. The sensitivity of this distance measure to scale means that vectors whose elements are further apart will have higher EDM values than vectors whose elements are close. In fact, the Euclidean distance satisfies the relation:

$$D_{Euc}(\alpha\vec{u},\alpha\vec{v}) = |\alpha|D_{Euc}(\vec{u},\vec{v}), \tag{11}$$

where $\alpha \in \mathbb{R}$ ($\alpha$ is any real number). We note that Eq. 11 does not mean that the EDM between two large vectors will be automatically be larger than the EMD between two relatively smaller vectors. The EDM is suitable for cases where the Euclidean distance between vectors is the dominant property of interest. Furthermore, the fact that the individual differences between corresponding elements in the two vectors are squared makes this similarity measure sensitive to outliers. A few outliers in any of the vectors can significantly shift the value of the index.

**2.4 City Block Distance Measure (CBDM)**

The city block distance (CBD) [45, 46] between two vectors $\vec{u}$ and $\vec{v}$ in an *n*-dimensional space is the sum of the absolute differences between the corresponding elements of the vectors. This distance measure is the $l_1$-norm of the difference of two vectors and it is given by:

$$D_{CB}(\vec{u},\vec{v}) = \|\vec{u}-\vec{v}\|_1 . \tag{12}$$

It may also be written as:

$$D_{CB}(\vec{u},\vec{v}) = \sum_{i=1}^{n} |u_i - v_i| . \tag{13}$$

The CBDM reflects an absolute comparison of the individual dimensions of the vectors being compared. Therefore, a pair of vectors whose corresponding elements are close to each other will yield a small city block distance value while a pair of vectors whose corresponding elements are far apart will yield a large city block distance value. Also, the CBDM is sensitive to scale such that vectors whose elements are further apart tend to yield high CBD values. The CBDM satisfies the relation:

$$D_{CB}(\alpha\vec{u},\alpha\vec{v}) = |\alpha| D_{CB}(\vec{u},\vec{v}) . \tag{14}$$

where $\alpha \in \mathbb{R}$. The CBDM is widely used in compressed sensing [47] and in assessing the differences in discrete probability distributions [48]. We note, in passing, that a generalized distance measure that encompasses both the EDM and the CBDM is the Minkowski distance measure (MDM) [49] given by:

$$D_{MD}(\vec{u},\vec{v}) = \left(\sum_{i=1}^{n} |u_i - v_i|^p\right)^{1/p} . \tag{15}$$

where *p* can be any integer from 1 to infinity (i.e., $p \in \mathbb{Z}^+$). The MDM reduces to the CBDM when $p=1$, and to the EDM when $p=2$. Also, note that the MDM is also sensitive to scale and satisfies the relation:

$$D_{MD}(\alpha\vec{u},\alpha\vec{v}) = |\alpha| D_{MD}(\vec{u},\vec{v}) . \tag{16}$$

where $D_{MD}(\vec{u},\vec{v})$ is the Minkowski distance measure between $\vec{u}$ and $\vec{v}$.

The relationship between the city block distance and the Euclidean distance is often given by the Cauchy-Schwarz inequality [50, 51]):

$$\|\vec{u}-\vec{v}\|_2 \leq \|\vec{u}-\vec{v}\|_1 \leq \sqrt{n}\|\vec{u}-\vec{v}\|_2, \tag{17}$$

where $\vec{u}, \vec{v} \in \mathbb{R}^n$. Equation (17) states that the Euclidean distance between any two vectors is less than or equal to the city block distance between them and the city block distance is less than or equal to the Euclidean distance multiplied by a factor of $\sqrt{n}$, $n$ being the dimension of the vectors. Also, another relation is often stated for the case in which all the elements of the vector $\vec{u}-\vec{v}$ are equal in value. This relation is stated as follows. For a vector $\vec{u}-\vec{v}$ with all elements equal,

$$\|\vec{u}-\vec{v}\|_1 = \sqrt{n}\|\vec{u}-\vec{v}\|_2. \tag{18}$$

Equation 18 states that for the special case in which the elements of the vector $\vec{u}-\vec{v}$ are all equal, the second inequality in Eq. 17 becomes an equality. We state here that the general equality relationship between the $l_1$- and $l_2$-norms is:

$$\|\vec{u}-\vec{v}\|_1^2 = \|\vec{u}-\vec{v}\|_2^2 + 2\sum_{i=1}^{n-1}\left\{(|u_i - v_i|)\left(\sum_{j=i+1}^{n}|u_j - v_j|\right)\right\}, \tag{19}$$

where $u_i$ and $v_i$ are the $i^{th}$ elements of vectors $\vec{u}$ and $\vec{v}$, respectively. To the best of our knowledge, Eq. 19 has not appeared in the literature prior to this article. Note that in Eq. 19, $\|\vec{u}-\vec{v}\|_1^2$ is the square of the $l_1$-norm of vector $(\vec{u}-\vec{v})$ and $\|\vec{u}-\vec{v}\|_2^2$ is the square of the $l_2$-norm of the vector. Thus, the CBD can be obtained by taking the square-root of the right-hand-side of Eq. 19. One way to evaluate the terms on the right-hand-side of Eq. 19 is to compute the outer product of the vector $|\vec{u}-\vec{v}|$ by itself to obtain a symmetric rank-one matrix $M$:

$$M = (|\vec{u}-\vec{v}|)|\vec{u}-\vec{v}|^T, \tag{20}$$

where $M$ is an $n \times n$ matrix, $n$ being the dimension of the vectors $\vec{u}$ and $\vec{v}$. Now, the first term on the right-hand-side of Eq. 19, $\|\vec{u}-\vec{v}\|_2^2$, is the trace of matrix $M$, and the last term without the factor of 2, i.e., $\sum_{i=1}^{n-1}\left\{(|u_i - v_i|)\left(\sum_{j=i+1}^{n}|u_j - v_j|\right)\right\}$, can be obtained by summing all the elements of matrix $M$ that fall above the main diagonal of the matrix. Although, the city block distance can be obtained from Eq. 19, computing this measure using Eq. 13 is a lot cheaper than computing it via Eq. 19. We present the proof of Eq. 19 in Appendix A and obtained Eq. 18 as a special case of Eq. 19.

**Hamming Distance Measure**

The Hamming distance [52, 53] measures the number of positions at which two equal-length vectors or two equal-length strings differ. For example, consider the strings: $s_1 = 2687\text{asdf}g1025$ and $s_2 = 2384\text{agdfm}1035$. The two strings are of equal lengths (i.e., each contains 13 characters).

Therefore, we can compute the Hamming distance between them. The number of corresponding positions at which the two strings differ is 5. Thus, the Hamming distance between the two strings is 5.

The Hamming distance is a distance metric in a Hamming space. This is because it fulfils all conditions associated with a metric [54]: (a) it satisfies the triangle inequality, (b) it is nonnegative (c) it is symmetric and (d) it equals zero if and only if the two strings being compared are identical. The Hamming distance finds widespread application in code theory, where it is used for error detection and error correction (channel coding).

### 3.3 Jaccard Similarity Index (JSI)

The Jaccard index [37, 55, 56, 57, 58] is primarily used for comparing sets and is defined as the ratio of the intersection of two sets to their union. The Jaccard similarity between any two sets A and B is:

$$S_{Jacc}(A,B) = \frac{|A \cap B|}{|A \cup B|}, \tag{21}$$

Since A and B are sets, each of them contains unique items. $A \cap B$ is the intersection of Sets A and B. It is a set that contains single record (avoiding duplication) of items that are present in both A and B. $A \cup B$ is the union of Sets A and B and it contains all items that are present in A or B or in both sets without duplication. $|A \cap B|$ is the number of items in $A \cap B$ and $|A \cup B|$ is the number of items in $A \cup B$. As an example, consider a bookstore with 100 different items and two customers A and B. The basket of Customer A contains pencil, crayon, ruler and notebook while the basket of Customer B contains pencil, pen and notebook. Then, $|A \cap B| = 2$ because only two items (pencil and notebook) are common to both customers. Also, $|A \cup B| = 5$ because there are five unique items in the baskets of both customers put together. Therefore, the JSI in this case is $2/5$. The JSI ranges from 0 to 1, with 0 indicating that there are no items common to the two sets; and 1 indicating that the two sets contain the same items (i.e., the two sets are equal). The JSI is often used in situations where only the presence or absence of features matters, such as in recommender systems.

### 3.4 Pearson Correlation Coefficient

Pearson correlation [6] measures the linear relationship between two vectors and is defined as:

$$S_{PCC}(\vec{u},\vec{v}) = \frac{\sum_{i=1}^{n}\left[(u_i - \bar{u})(v_i - \bar{v})\right]}{\sqrt{\left[\sum_{i=1}^{n}(u_i - \bar{u})^2\right]\left[\sum_{i=1}^{n}(u_i - \bar{v})^2\right]}}, \tag{22}$$

where $\bar{u}$ and $\bar{v}$ are the mean values of the vectors $\vec{u}$ and $\vec{v}$, respectively. The result ranges from $-1$ (perfect negative correlation) to 1 (perfect positive correlation), with 0 indicating no linear relationship.

## 4. The Joint Distance Measure (JDM)

In this article, we introduce a vector similarity measure known as the joint distance measure (JDM). This measure is formed by the combination of the cosine similarity measure and any of the Minkowski distance metric. The JDM is defined as:

$$D_{JDM}(\vec{u},\vec{v}) = \left(\sum_{i=1}^{n} |u_i - v_i|^p\right)^{1/p} |\cos\theta - 2|, \tag{23}$$

where $\theta$ is the angle between the vectors $\vec{u}$ and $\vec{v}$. The value of $p$ dictates the distance metric that is combined with the cosine similarity measure. When $p=1$, the city block distance metric is combined with the cosine similarity measure and when $p=2$, the Euclidean distance metric is combined with the cosine similarity measure. Equation 23 can also be expressed as:

$$D_{JDM}(\vec{u},\vec{v}) = D_{MD} |S_{\cos} - 2|, \tag{24}$$

where $D_{MD}$ is the Minkowski distance metric and $S_{\cos}$ is the cosine similarity measure. When $p=2$, the JDM becomes:

$$D_{JDM}(\vec{u},\vec{v}) = D_{Euc} |S_{\cos} - 2|. \tag{25}$$

The values of JDM is such that $D_{MD} \leq D_{JDM} \leq 3D_{MD}$. That is, for any pair of vectors, the minimum achievable value of the JDM is the value of the corresponding MDM and the maximum attainable value of the JDM is thrice the value of the corresponding MDM. The minimum value is obtained only when the pair of vectors are parallel and pointing in the same direction. The maximum value is obtained when the vectors are parallel but pointing in opposite directions. When the vectors are orthogonal, then the JDM is twice the size of the corresponding MDM (i.e., $D_{JDM} = 2D_{MD}$). Note that the joint distance measure may be posed as a similarity measure as:

$$S_{JSM} = \frac{1}{1 + D_{JDM}}, \tag{26}$$

where $S_{JSM}$ is the joint similarity measure (JSM). Similar equation for the Euclidean distance can be found in [37]. The range of the JSM is $0 < S_{JDM} \leq 1$, with 1 indicating completely similar vectors. The value of 1 can be attained if and only if the vectors being compared are identical. Values of JSM closer to zero indicate that the vectors are highly dissimilar with large spatial and angular distances between the vectors.

Although, the Minkowski distance measure is a distance metric, the cosine similarity measure and the joint distance measure (JDM) are not. In Appendix B, we show that the JDM is not a metric.

**Cost Implications of the Joint Distance Measure**

We note that computing the JDM requires the computation of the desired distance metric and the cosine similarity measure. However, for the case of the Euclidean distance metric ($p = 2$), we can obtain the Euclidean distance measure from the components of the cosine similarity measure. We note that the square of the Euclidean distance measure can be written as [50]:

$$\|\vec{u} - \vec{v}\|_2^2 = \|\vec{u}\|_2^2 - 2\vec{u} \cdot \vec{v} + \|\vec{v}\|_2^2. \qquad (27)$$

However, from Eq. 2, we observe that computing the CSM requires the computation of $\|\vec{u}\|$, $\|\vec{v}\|$ and $\vec{u} \cdot \vec{v}$, which are also the components of Eq. 26. Once $\|\vec{u}\|$, $\|\vec{v}\|$ and $\vec{u} \cdot \vec{v}$ are obtained when computing the CSM, they can be used in Eq. 26 to compute the square of the EDM. Thus, in this case, instead of computing the EDM from Eq. 10, we may compute it by taking the square-root of Eq. 26 since the terms in Eq. 26 are already available. If the city block distance measure (with $p = 1$) is used in place of the EDM to compute the JDM, we suggest that Eq. 13 should be used to compute the CBDM because it is relatively cheaper to compute the sum of absolute values.

**Examples**

In this section, we present examples that illustrate the behavior of the JDM and its comparison to the CSM, EDM, CBDM and the DPSM.

Example 1 – Illustrating the behaviors of some similarity measures

In this example, we considered the vectors $\vec{r} = [-1, 0]$, $\vec{s} = [1, 0]$, $\vec{t} = [4, 0]$, $\vec{u} = [8, 0]$ $\vec{v} = [10, 0]$, $\vec{w} = [0, 2]$, $\vec{x} = [-1, 5]$, $\vec{y} = [2, -4]$ and $\vec{z} = [1, 1]$; and computed the similarity/distance measures between some pairs of these vectors as shown in Table 1. The measures considered are the CSM, DPSM, CBDM, EDM, JDM 1 (JDM with $p = 1$) and JDM 2 (JDM with $p = 2$). The exercise was performed to illustrate the strengths and weaknesses of these measures.

| Table 1 – Distance/Similarity Measures illustrating the behaviors of the different measures | | | | | | | | |
|---|---|---|---|---|---|---|---|---|
| Case | Vector pair | $\theta$ | CSM | DPSM | CBDM | EDM | JDM 1 | JDM 2 |
| 1 | $\vec{r}$ and $\vec{s}$ | 180º | –1 | –1 | 2 | 2 | 6 | 6 |
| 2 | $\vec{s}$ and $\vec{t}$ | 0º | 1 | 4 | 3 | 3 | 3 | 3 |
| 3 | $\vec{u}$ and $\vec{v}$ | 0º | 1 | 80 | 2 | 2 | 2 | 2 |
| 4 | $\vec{s}$ and $\vec{w}$ | 90º | 0 | 0 | 3 | 2.236 | 6 | 4.472 |
| 5 | $\vec{t}$ and $\vec{w}$ | 90º | 0 | 0 | 6 | 4.472 | 12 | 8.944 |

| 6 | $\vec{s}$ and $\vec{z}$ | 45° | 0.707 | 1 | 1 | 1 | 1.293 | 1.293 |
| 7 | $\vec{x}$ and $\vec{y}$ | 164.8° | –0.965 | –22 | 12 | 9.487 | 35.577 | 28.126 |

Table 1 displays seven cases used to illustrate the behaviors of the similarity/distance measures. In Case 1, vectors $\vec{r}$ and $\vec{s}$ are parallel and pointing in opposite directions and the JDM 1 and JDM 2 are three times larger than the CBDM and EDM, respectively. This supports the notion that when vectors are pointing in opposite direction (angular difference is 180°), the JDM yields values that are three times larger than the respective Minkowski distance measures. In Cases 2 and 3, in which the vectors in each pair are parallel to each other and pointing in the same direction, the JDM 1 and JDM 2 are the same as the CBDM and EDM, respectively. These two cases support the fact that for parallel vectors pointing in the same direction (angular difference is zero), the joint distance measures yield the same values as the values from the corresponding Minkowski distance measures. In Cases 4 and 5, the vectors in each pair are at 90° to each other. In these cases, the JDM yields values that are twice as large as the values from the corresponding Minkowski distance measure. Thus, in summary, the JDM gives a multiple of the Minkowski distance measure, with the smallest multiple of 1 obtained when the vectors in the pair are parallel and pointing in the same direction. This multiplicative factor gradually increases to 2 as the angle between the pair of vectors increase from 0 to 90°. As the angular distance between the pair of vectors increases from 90° to 180°, the multiplicative factor increases from 2 to 3.

We observe that the values of the Minkowski distance measures in Case 1 are the same as those in Case 3. However, the values of the JDM in these two cases are different due to differences in angular separation of the pairs of vectors in the two cases. In Cases 4 and 5, the angular distances are the same so that the CSM resulted in the same value in both cases. However, the values from the JDMs (JDM 1 and JDM 2) are different in these two cases. These examples show that the joint distance measures have more distinguishing features than the CSM and the MDM. In fact, if two vectors are at the same angular distance from a third vector (call this the index vector) but their spatial distances to the index vector are different, the vectors will have different JDM values with respect to the index vector. Also, if two vectors have the same spatial distance but different angular distances to an index vector, the vectors will have different JDM values to the index vector. Case 6 is one in which the angular distance between the pair of vectors is between 0 and 90. In this case, the multiplicative factor is between 1 and 2.

Example 2 – Illustrating the inconsistency of the dot product similarity measure

This example is used to illustrate the inconsistent behavior of the dot product similarity measure (DPSM). We consider two scenarios under which the similarities between an index vector $\vec{z}_0$ and four other vectors are computed. The four vectors remain the same under both scenarios. However, the vector used as the index vector in first scenario is different from that used as the index vector in the second scenario.

In Scenario 1, the index vector is $\vec{z}_0 = [-1, 2]$, while other vectors are $\vec{z}_1 = [-2, 3]$, $\vec{z}_2 = [-12, 16]$, $\vec{z}_3 = [-5, -8]$ and $\vec{z}_4 = [12, 16]$. The DPSM and other measures were computed and displayed in Table 2.

| Table 2 – Distance/Similarity Measures illustrating the behaviors of the different measures (Scenario 1) | | | | | | | | |
|---|---|---|---|---|---|---|---|---|
| Case | Vector pair | $\theta$ | CSM | DPSM | CBDM | EDM | JDM 1 | JDM 2 |
| 1 | $\vec{z}_0$ and $\vec{z}_0$ | 0º | 1 | 5 | 0 | 0 | 0 | 0 |
| 2 | $\vec{z}_0$ and $\vec{z}_1$ | 7.12º | 0.992 | 8 | 2 | 1.414 | 2.015 | 1.425 |
| 3 | $\vec{z}_0$ and $\vec{z}_2$ | 10.3º | 0.984 | 44 | 25 | 17.8 | 25.4 | 18.09 |
| 4 | $\vec{z}_0$ and $\vec{z}_3$ | 121.4º | –0.521 | –11 | 14 | 10.77 | 35.3 | 27.16 |
| 5 | $\vec{z}_0$ and $\vec{z}_4$ | 63.4º | 0.447 | 20 | 27 | 19.11 | 41.93 | 29.67 |

Five cases were considered and shown in Table 2. Case 1, in which the similarity of the index vector $\vec{z}_0$ to itself was evaluated, was included as a base-case to gauge the other cases. As shown on the table, we got the largest DPSM in Case 3 followed by Case 5. The lowest value of DPSM was obtained in Case 3. However, merely looking at the set of vectors, we can tell that the closest amongst the four other vectors to the index vector $\vec{z}_0$ is the vector $\vec{z}_1$. Therefore, DPSM erroneously implies that $\vec{z}_2$ and $\vec{z}_4$ are closer to $\vec{z}_0$ than $\vec{z}_1$ is. In fact, if we try to determine how similar $\vec{z}_0$ is to itself (Case 1) using the DPSM, we get a similarity value of 5, indicating that $\vec{z}_0$ is more similar to vectors $\vec{z}_1$, $\vec{z}_2$ and $\vec{z}_4$, than it is to itself. This result is erroneous, indicating that the dot product is not a reliable similarity measure. We note, in addition, that any vector that is larger in size than the index vector and in the same direction as the index vector will be evaluated as being more similar to the index vector than the index vector is to itself. This is because the DPSM is erroneously using the size of the vector being compared to the index vector as a major determiner of similarity. In general, there are an infinite number of vectors that are in the same direction as the index vector and at the same time larger in size than the index vector.

Table 3 shows the ranking of the closeness of the vectors to the index vector based on the various measures. Clearly, the CBDM and the EDM have correctly ranked all the vectors based on the spatial distance only while the CSM have also correctly ranked the vectors based on the angular distance only. However, the correct ranking based on joint spatial and angular distances was given only by JDM 1 and JDM 2.

| Table 3 – Ranking of vectors based on similarity to index vector (Scenario 1) | | | | | | | | |
|---|---|---|---|---|---|---|---|---|
| Case | Vector pair | $\theta$ | CSM | DPSM | CBDM | EDM | JDM 1 | JDM 2 |
| 1 | $\vec{z}_0$ and $\vec{z}_0$ | 0º | 1 | 4 | 1 | 1 | 1 | 1 |
| 2 | $\vec{z}_0$ and $\vec{z}_1$ | 7.12º | 2 | 3 | 2 | 2 | 2 | 2 |
| 3 | $\vec{z}_0$ and $\vec{z}_2$ | 10.3º | 3 | 1 | 4 | 4 | 3 | 3 |
| 4 | $\vec{z}_0$ and $\vec{z}_3$ | 121.4º | 5 | 5 | 3 | 3 | 4 | 4 |
| 5 | $\vec{z}_0$ and $\vec{z}_4$ | 63.4º | 4 | 2 | 5 | 5 | 5 | 5 |

In Scenario 2, the index vector was $\vec{z}_0 = [-10, 15]$, while other vectors were $\vec{z}_1 = [-2, 3]$, $\vec{z}_2 = [-12, 16]$, $\vec{z}_3 = [-5, -8]$ and $\vec{z}_4 = [12, 16]$. The DPSM and other measures were computed and displayed in Table 4.

| Table 4 – Distance/Similarity Measures illustrating the behaviors of the different measures (Scenario 2) | | | | | | | | |
|---|---|---|---|---|---|---|---|---|
| Case | Vector pair | $\theta$ | CSM | DPSM | CBDM | EDM | JDM 1 | JDM 2 |
| 1 | $\vec{z}_0$ and $\vec{z}_0$ | 0° | 1 | 325 | 0 | 0 | 0 | 0 |
| 2 | $\vec{z}_0$ and $\vec{z}_1$ | 0° | 1 | 65 | 20 | 14.42 | 20 | 14.42 |
| 3 | $\vec{z}_0$ and $\vec{z}_2$ | 3.18° | 0.998 | 360 | 3 | 2.236 | 3.005 | 2.24 |
| 4 | $\vec{z}_0$ and $\vec{z}_3$ | 114.3° | –0.412 | –70 | 28 | 23.54 | 67.52 | 56.76 |
| 5 | $\vec{z}_0$ and $\vec{z}_4$ | 70.6° | 0.333 | 120 | 23 | 22.02 | 38.35 | 36.72 |

The table shows that the largest DPSM value was obtained in Case 3 while the lowest was obtained in Case 4. Also, looking at the set of vectors, we can see that among the four other vectors, the vector $\vec{z}_2$ (Case 3) is the most similar vector to the index vector $\vec{z}_0$. In this example, the DPSM indicated correctly that the vector $\vec{z}_2$ is more similar to the index vector $\vec{z}_0$ than any of the three other vectors is to the index vector. The reason for this is that, in this scenario, vector $\vec{z}_2$ is large relative to the other vectors and it is almost parallel to the index vector $\vec{z}_0$, which is also large relative to the vectors in the group. Thus, the dot product of $\vec{z}_0$ with $\vec{z}_2$ is larger than the dot product of $\vec{z}_0$ with each of the other three vectors. This is however, not an indication of similarity. It is more an indication of the relative size of $\vec{z}_2$. Any vector that is larger in size than $\vec{z}_2$ and parallel or almost parallel to $\vec{z}_0$ will have a bigger dot product with $\vec{z}_0$ than the dot product of $\vec{z}_0$ with $\vec{z}_2$. Also, a closer look at Table 4 shows that the DPSM erroneously indicated that $\vec{z}_2$ is more similar to $\vec{z}_0$ than $\vec{z}_0$ is to itself. The ranking of the vectors based on the measures are displayed in Table 5.

| Table 5 – Ranking of vectors based on similarity to index vector (Scenario 2) | | | | | | | | |
|---|---|---|---|---|---|---|---|---|
| Case | Vector pair | $\theta$ | CSM | DPSM | CBDM | EDM | JDM 1 | JDM 2 |
| 1 | $\vec{z}_0$ and $\vec{z}_0$ | 0° | 1 | 2 | 1 | 1 | 1 | 1 |
| 2 | $\vec{z}_0$ and $\vec{z}_1$ | 0° | 1 | 4 | 3 | 3 | 3 | 3 |
| 3 | $\vec{z}_0$ and $\vec{z}_2$ | 3.18° | 3 | 1 | 2 | 2 | 2 | 2 |
| 4 | $\vec{z}_0$ and $\vec{z}_3$ | 114.3° | 5 | 5 | 5 | 5 | 5 | 5 |
| 5 | $\vec{z}_0$ and $\vec{z}_4$ | 70.6° | 4 | 3 | 4 | 4 | 4 | 4 |

# 6. Conclusion

Vector similarity measures are foundational to numerous algorithms and applications across machine learning, data mining, NLP, and information retrieval. Each measure has its strengths and weaknesses, and the choice of measure depends on the specific task and the nature of the data. In this article, we proposed the Joint Distance Measure (JDM) that is combines the advantages of the Minkowski distance measure (MDM) and the cosine similarity measure (CSM). We showed the strengths of the JDM as compared with other similar similarity measures.

---

**References**

1. James Abello, Panos M. Pardalos, and Mauricio. G. C. Resende. Handbook of Massive Data Sets. *Kluwer Academic Publishers,* USA, 2002.
2. Witold Abramowicz. Knowledge-Based Information Retrieval and Filtering from the Web. *Springer,* Boston, MA, 2003. doi:10.1007/978-1-4757-3739-4.
3. Hyung Jun Ahn. A new similarity measure for collaborative filtering to alleviate the new user cold-starting problem. *Information Sciences* 178 (1), 37–51, 2008. https://doi.org/10.1016/j.ins.2007.07.024.
4. Gediminas Adomavicius, Bamshad Mobasher, Francesco Ricci, and Alexander Tuzhilin. Context-aware recommender systems. *AI Magazine* 32 (3), 67–80, 2011. https://doi.org/10.1609/ aimag.v32i3.2364.
5. Charu C. Aggarwal. Recommender Systems: The Textbook. *Springer Publishing Company Incorporated*, 2016.
6. Karl Pearson. Note on Regression and Inheritance in the Case of Two Parents, *Proceedings of the Royal Society of London* Series I (58), pp. 240-242, 1895.
7. Hiyan K. Alshawi. The core language engine. *MIT press*, 1992.
8. H. Ahonen, O. Heinonen, M. Klemettinen and A. I. Verkamo. Applying data mining techniques for descriptive phrase extraction in digital document collections, In *Proceedings IEEE International Forum on Research and Technology Advances in Digital Libraries -ADL'98-*, Santa Barbara, CA, USA, 1998, pp. 2-11, doi: 10.1109/ADL.1998.670374.
9. Shiliang Sun, Chen Luo, Junyu Chen. A review of natural language processing techniques for opinion mining systems, *Inf. Fusion*, vol. 36, pp. 10-25, 2017.
10. Prashant Gupta, Aman Goswami, Sahil Koul, Kashinath Sartape. IQS-intelligent querying system using natural language processing, *Proc. Int. Conf. Electron. Commun. Aerosp. Technol. ICECA 2017*, vol. 2017–Jan., pp. 410-413, 2017.
11. Basemah Alshemali, Jugal Kalita. Improving the reliability of deep neural networks in NLP: A review. *Knowl-Based Syst* 191:105210, 2020.
12. Tomas Mikolov, Kai Chen, Greg Corrado, and Jeffrey Dean. Efficient Estimation of Word Representations in Vector Space. *arXiv*:1301.3781, 2013.

13. Tomas Mikolov, Ilya Sutskever, Kai Chen, Greg Corrado, and Jeffrey Dean. Distributed representations of words and phrases and their compositionality. *Advances in Neural Information Processing Systems* - Volume 2, Pages 3111 - 3119, 2013.
14. Soubraylu Sivakumar, Lakshmi Sarvani Videla, T Rajesh Kumar, J. Nagaraj, Shilpa Itnal, and D. Haritha. Review on Word2Vec Word Embedding Neural Net, *2020 International Conference on Smart Electronics and Communication*, Trichy, India, pp. 282-290, 2020. doi: 10.1109/ICOSEC49089.2020.9215319
15. Jeffrey Pennington, Richard Socher, and Christopher Manning. GloVe: Global Vectors for Word Representation. *Proceedings of the 2014 Conference on Empirical Methods in Natural Language Processing*, Doha, Qatar: Association for Computational Linguistics: 1532–1543, 2014. doi:10.3115/v1/D14-1162.
16. Jacob Devlin, Ming-Wei Chang, Kenton Lee, and Kristina Toutanova. BERT: Pre-training of Deep Bidirectional Transformers for Language Understanding. *arXiv*:1810.04805v2, 2018.
17. Kawin Ethayarajh. How Contextual are Contextualized Word Representations? Comparing the Geometry of BERT, ELMo, and GPT-2 Embeddings, *arXiv*:1909.00512, 2019.
18. Anna Rogers, Olga Kovaleva, and Anna Rumshisky. A Primer in BERTology: What We Know About How BERT Works. *Transactions of the Association for Computational Linguistics*, 8: 842–866. arXiv:2002.12327,
2020. doi:10.1162/tacl_a_00349. S2CID 211532403.
19. William B. Frakes, Ricardo Baeza-Yates. Information Retrieval Data Structures & Algorithms. *Prentice-Hall, Inc*. ISBN 978-0-13-463837-9, 1992.
20. Jonathan Foote. An overview of audio information retrieval. *Multimedia Systems*, 7: *2*–10, 1999. doi:10.1007/s005300050106.
21. Abby A. Goodrum. Image Information Retrieval: An Overview of Current Research. *Informing Science*, Vol 3(2): 63-66, 2000.
22. Amit Singhal. Modern Information Retrieval: A Brief Overview. *IEEE Data Eng. Bull.*, 24 (4): *35*–43, 2001.
23. Jöran Beel, Bela Gipp, and Jan-Olaf Stiller. Information retrieval on mind maps - What could it be good for? *Proceedings of the 5th International Conference on Collaborative Computing: Networking, Applications and Worksharing*, Washington, DC, 2009.
24. Rui Xu, and D. Wunsch, Survey of clustering algorithms. *IEEE Transactions on Neural Networks*, 16(3), 645–678, 2005. https://doi.org/10.1109/TNN.2005.845141
25. Anil K. Jain. Data clustering: 50 years beyond K-means. *Pattern Recognition Letters*, Volume 31, Issue 8, Pages 651-666, 2010.
26. Dongkuan Xu, and Yingjie Tian, A comprehensive survey of clustering algorithms. *Annals of Data Science*, Vol 2, 165–193, 2015. https://doi.org/10.1007/s40745-015-0040-1
27. Hui Yin, Amir Aryani, Stephen Petrie, Aishwarya Nambissan, Aland Astudillo, and Shengyuan Cao. A Rapid Review of Clustering Algorithms. *arxiv*.org/pdf/2401.07389, 2024.
28. Abiodun M. Ikotun, Absalom E. Ezugwu, Laith Abualigah, Belal Abuhaija, and Jia Heming. K-means clustering algorithms: A comprehensive review, variants analysis, and advances in the era of big data. *Information Sciences*, Volume 622, Pages 178-210, 2023. https://doi.org/10.1016/j.ins.2022.11.139

29. Xingcheng Ran, Yue Xi, Yonggang Lu, Xiangwen Wang, and Zhenyu Lu. Comprehensive survey on hierarchical clustering algorithms and the recent developments. *Artif Intell Rev* 56, 8219–8264, 2023. https://doi.org/10.1007/s10462-022-10366-3.

30. Xudong Jiang, and Alvin H. Wah, Constructing and training feed-forward neural networks for pattern classification. *Pattern Recognition*, Volume 36, Issue 4, April 2003, Pages 853-867, 2003.

31. Guobin Ou, Yi Lu Murphey. Multi-class pattern classification using neural networks. *Pattern Recognition*, Volume 40, Issue 1, Pages 4-18, 2007.

32. Zhaodong Wu, Jun Zhang, and Shengliang Hu. Review on Classification Algorithm and Evaluation System of Machine Learning. *13th International Conference on Intelligent Computation Technology and Automation*, Xi'an, China, pp. 214-218, 2020. doi: 10.1109/ICICTA51737.2020.00052.

33. Rahul Sharma, Raphael Hussung, Andreas Keil, Fabian Friederich, Thomas Fromenteze, and Mohsen Khalily. Performance Analysis of Classification Algorithms for Millimeter-wave Imaging. *16th European Conference on Antennas and Propagation*, Madrid, Spain, pp. 1-5, 2022. doi: 10.23919/EuCAP53622.2022.9769429.

34. Jiawei Han, Micheline Kamber, and Jian Pei. Chapter 2 - Getting to Know Your Data, *Data Mining* (Third Edition), The Morgan Kaufmann Series in Data Management Systems, Pages 39-82, 2012.

35. Saprativa Bhattacharjee, Anirban Das, Ujjwal Bhattacharya, Swapan K. Parui, and Sudipta Roy. Sentiment analysis using cosine similarity measure. *2015 IEEE 2nd International Conference on Recent Trends in Information Systems*, Kolkata, India, pp. 27-32, 2015. doi: 10.1109/ReTIS.2015.7232847.

36. Fethi Fkih. Similarity measures for Collaborative Filtering-based Recommender Systems: Review and experimental comparison. *Journal of King Saud University - Computer and Information Sciences*, Volume 34(9), Pages 7645-7669, 2022.

37. Pradipto Chowdhury, and Bam Bahadur Sinha. Evaluating the Effectiveness of Collaborative Filtering Similarity Measures: A Comprehensive Review. *International Conference on Machine Learning and Data Engineering*, 2023.

38. Mohammad Yasser, Khaled F. Hussain, and Samia A. Ali. Comparative Analysis of Similarity Methods in High-Dimensional Vectors: A Review. *International Conference on Artificial Intelligence Science and Applications in Industry and Society*, Galala, Egypt, pp. 1-6, 2023. doi: 10.1109/CAISAIS59399.2023.10270776.

39. Brendan O'Connor. Cosine similarity, Pearson correlation, and OLS coefficients, *AI and Social Science*, Brendan O'Connor, Posted on March 13, 2012. https://brenocon.com/blog/2012/03/cosine-similarity-pearson-correlation-and-ols-coefficients/. Retrieved 20-Mar-2025.

40. Frederik vom Lehn. Understanding Vector Similarity for Machine Learning, Published in *Advanced Deep Learning*, Posted on Oct 7, 2023, medium.com 2023. https://medium.com/advanced-deep-learning/understanding-vector-similarity-b9c10f7506de. Retrieved 20-Mar-2025.

41. Oracle. Dot Product Similarity, *Oracle AI Vector Search User's Guide*, 2023 Release: https://docs.oracle.com/en/database/oracle/oracle-database/23/vecse/dot-product-similarity.html. Retrieved 20-Mar-2025
42. Roie Schwaber-Cohen. Vector Similarity Explained. *Pinecone*, Posted on June 30, 2023. https://www.pinecone.io/learn/vector-similarity/, Retrieved 20-Mar-2025
43. Robert J. Bell. An Elementary Treatise on Coordinate Geometry of Three Dimensions (2nd ed.), *Macmillan*, pp. 57–61, 1914.
44. John Tabak. Geometry: The Language of Space and Form. *Facts on File math library, Infobase Publishing*, p. 150, 2014. ISBN 978-0-8160-6876-0.
45. Erdogan. S. Şuhubi. Chapter V: Metric Spaces, Functional Analysis. *Springer Netherlands*, pp. 261–356, 2003. doi:10.1007/978-94-017-0141-9_5
46. Pavel Zezula, Giuseppe Amato, Vlastislav Dohnal, and Michal Batko. Similarity Search: The Metric Space Approach. *Advances in Database Systems*, Springer, p. 10, 2006. doi:10.1007/0-387-29151-2.
47. David L. Donoho. For most large underdetermined systems of linear equations the minimal-norm solution is also the sparsest solution. *Communications on Pure and Applied Mathematics*. 59 (6): 797–829, 2006. doi:10.1002/cpa.20132. S2CID 8510060.
48. Kian Huat Lim, Luciana Ferraris, Madeleine E. Filloux, Benjamin J. Raphael, and William G. Fairbrother. Using positional distribution to identify splicing elements and predict pre-mRNA processing defects in human genes. *Proceedings of the National Academy of Sciences of the United States of America*, 108 (27): 11093–11098, 2011. doi:10.1073/pnas.1101135108.
49. Hermann Minkowski. Geometrie der Zahlen (in German). *Leipzig and Berlin: R. G. Teubner*. JFM 41.0239.03. MR 0249269, 1910.
50. Gene H. Golub, and Charles F. Van Loan. Matrix Computations (Third ed.). Baltimore: *The Johns Hopkins University Press*. p. 53, 1996. ISBN 0-8018-5413-X.
51. J. Micheal Steele. The Cauchy–Schwarz Master Class: An Introduction to the Art of Mathematical Inequalities. *Cambridge University Press*, 2004.
52. R. W. Hamming. Error detecting and error correcting codes. *The Bell System Technical Journal*, 29 (2): 147–160,1950. doi:10.1002/j.1538-7305.
53. Ayman Jarrous, and Benny Pinkas. Secure Hamming distance based computation and its applications. *Applied Cryptography and Network Security, Lecture Notes in Computer Science*, Vol. 5536. Berlin, Heidelberg: Springer. pp. 107–124, 2009. doi:10.1007/978-3-642-01957-9.
54. Derek J. Robinson. An Introduction to Abstract Algebra. *Walter de Gruyter*, pp. 255–257, 2003.
55. Paul Jaccard. Étude comparative de la distribution florale dans une portion des Alpes et des *Jura. Bulletin de la Société vaudoise des sciences naturelles* (in French). 37 (142): *547*–579, 1901.
56. Paul Jaccard. The Distribution of the Flora in the Alpine Zone.1, *New Phytologist*, 11 (2): *37*–50, 1912. doi:10.1111/j.1469-8137.1912.tb05611.x.
57. Ryan Moulton, and Yunjiang Jiang. Maximally Consistent Sampling and the Jaccard Index of Probability Distributions. *2018 IEEE International Conference on Data Mining*. pp. 347–356. arXiv:1809.04052, 2018. doi:10.1109/ICDM.2018.00050.

58. Sven Kosub. A note on the triangle inequality for the Jaccard distance. *Pattern Recognition Letters*. 120: 36 – 38. arXiv:1612.02696, 2019. doi:10.1016/j.patrec.2018.12.007.

---

**Appendix A – Proof of the relation between the $l_1$-norm and the $l_2$-norm**

In this appendix, we prove the result presented in Eq. 19. It is important to note here that we have searched the open literature for a general equality relation between the $l_1$-norm and the $l_2$-norm, but could not find any other than that shown in Eq. 18. It is however necessary to present the basis for which we claim Eq. 19 to be a valid. While we present the relation in Eq. 19 and its proof in this paper, it may well be that these are a subset of a broader theorem and proof that have previously been established in Geometry or in the theory of normed spaces. But until we see that, we present our proof of the relation here.

Consider two real numbers $a_1$ and $a_2$ (i.e., $a_1, a_2 \in \mathbb{R}$). The square of the sum of these numbers is

$$(a_1 + a_2)^2 = a_1^2 + a_2^2 + 2a_1 a_2. \tag{A1}$$

We may repeat this for three real numbers $a_1$, $a_2$ and $a_3$ as

$$(a_1 + a_2 + a_3)^2 = a_1^2 + a_2^2 + a_3^2 + 2(a_1 a_2 + a_1 a_3 + a_2 a_3). \tag{A2}$$

Now, for $n$ real numbers $a_1$, $a_2$, ..., $a_n$, we have

$$(a_1 + a_2 + \cdots + a_n)^2 = a_1^2 + a_2^2 + a_3^2 + \cdots + a_n^2 + 2 \begin{pmatrix} a_1 a_2 + a_1 a_3 + a_1 a_4 + \cdots + a_1 a_n + \\ a_2 a_3 + a_2 a_4 + a_2 a_5 + \cdots + a_2 a_n + \\ a_3 a_4 + a_3 a_5 + a_3 a_6 + \cdots + a_3 a_n + \\ \cdots + a_{n-1} a_n \end{pmatrix} \tag{A3}$$

Equation A3 may be written as

$$(a_1 + a_2 + \cdots + a_n)^2 = a_1^2 + a_2^2 + a_3^2 + \cdots + a_n^2 + 2 \begin{bmatrix} a_1 (a_2 + a_3 + a_4 + \cdots + a_n) + \\ a_2 (a_3 + a_4 + a_5 + \cdots + a_n) + \\ a_3 (a_4 + a_5 + a_6 + \cdots + a_n) + \\ \cdots + a_{n-1} a_n \end{bmatrix}. \tag{A4}$$

In summation notation, we may express Eq. A4 as

$$\left(\sum_{i=1}^{n} a_i\right)^2 = \sum_{i=1}^{n} a_i^2 + 2\sum_{i=1}^{n-1}\left\{(a_i)\left(\sum_{j=i+1}^{n} a_j\right)\right\}. \tag{A5}$$

We emphasize that Eq. A5 is applicable to any set of real numbers which may be all positive, all negative or a mix of positive and negative numbers. If we take the real numbers $a_1$, $a_2$, …, $a_n$ as elements of an $n$-dimensional vector $\vec{a}$, then $\sum_{i=1}^{n} a_i^2 = \|\vec{a}\|_2^2$ but $\left(\sum_{i=1}^{n} a_i\right)^2$ may not necessary be equal to $\|\vec{a}\|_1^2$ because the set of real numbers $a_1$, $a_2$, …, $a_n$ may contain a mix of positive and negative numbers. Now, if we restrict our attention to only positive numbers and consider that the set of numbers $|a_1|$, $|a_2|$, …, $|a_n|$ are the elements of vector $\vec{a}$, then we may repeat the steps in Eqs. A1 to A5 to obtain

$$\left(\sum_{i=1}^{n} |a_i|\right)^2 = \sum_{i=1}^{n} a_i^2 + 2\sum_{i=1}^{n-1}\left\{(|a_i|)\left(\sum_{j=i+1}^{n} |a_j|\right)\right\}. \tag{A6}$$

Equation A6 is a special case of Eq. A5 and in this case, $\left(\sum_{i=1}^{n} |a_i|\right)^2 = \|\vec{a}\|_1^2$. Thus, Eq. A6 can be rewritten as

$$\|\vec{a}\|_1^2 = \|\vec{a}\|_2^2 + 2\sum_{i=1}^{n-1}\left\{(|a_i|)\left(\sum_{j=i+1}^{n} |a_j|\right)\right\}. \tag{A7}$$

Undoubtedly, Eq. A7 establishes a relationship between the $l_1$ and $l_2$-norms of a vector $\vec{a}$ whose elements are absolute values of some real numbers. If we take $\vec{a}$ to be the absolute value of the difference between the two vectors $\vec{u}$ and $\vec{v}$ so that $\vec{a} = |\vec{u} - \vec{v}|$ and $a_i = |u_i - v_i|$ for every index $i \in \{1, 2, \cdots, n\}$, then Eq. A7 may be written as

$$\|\vec{u} - \vec{v}\|_1^2 = \|\vec{u} - \vec{v}\|_2^2 + 2\sum_{i=1}^{n-1}\left\{(|u_i - v_i|)\left(\sum_{j=i+1}^{n} |u_j - v_j|\right)\right\}, \tag{A8}$$

which is the same as Eq. 19.

If in Eq. A7, all the elements in vector $\vec{a}$ are equal in value (i.e., $a_1 = a_2 = \cdots = a_n = a$), then the summation term (last term excluding the factor of 2) in Eq. A7 reduces to $\dfrac{n(n-1)a^2}{2}$ and the $l_2$-norm reduces to $na^2$ so that the equation becomes

$$\|\vec{a}\|_1^2 = na^2 + 2\frac{n(n-1)a^2}{2}. \tag{A9}$$

Further simplification of Eq. A9 yields:

$$\|\vec{a}\|_1^2 = n^2 a^2. \tag{A10}$$

Since for this special case of a vector with all elements equal, $\|\vec{a}\|_2^2 = na^2$, we may rewrite Eq. A10 as:

$$\|\vec{a}\|_1^2 = n\|\vec{a}\|_2^2. \tag{A11}$$

Equation 18 in this article can be obtained from by taking the square-root of both sides of Eq. A11.

We test our derivations using the following vectors $\vec{a} = [-2,\ 4,\ 1,\ -8]$, $\vec{u} = [-7,\ -1,\ 3,\ -2]$ and $\vec{v} = [-5,\ 2,\ -1,\ -4]$. The first vector will be used to test the more general relation given in Eq. A5 while the other two vectors will be used to test the relation presented in Eq. A8. The first term in Eq. A5 is evaluated as

$$\left(\sum_{i=1}^{n} a_i\right)^2 = (-2+4+1-8)^2 = (-5)^2 = 25.$$

The second term in Eq. A5 is evaluated as:

$$\sum_{i=1}^{n} a_i^2 = (-2)^2 + 4^2 + 1^2 + (-8)^2 = 4+16+1+64 = 85.$$

The last term in Eq. A5 is evaluated as

$$2\sum_{i=1}^{n-1}\left\{(a_i)\left(\sum_{j=i+1}^{n} a_j\right)\right\} = 2\left[-2(4+1-8)+4(1-8)+1(-8)\right]$$
$$= 2\left[-2(-3)+4(-7)+1(-8)\right]$$
$$= 2(6-28-8) = 2(-30) = -60$$

Since 25 is the same as $85-60$, we note that the left-hand-side of Eq. A5 evaluates to the right-hand-side of the equation.

Now, to evaluate the terms in Eq. A8, we first compute $\vec{u} - \vec{v} = [-2,\ -3,\ 4,\ 2]$. Then, we evaluate the first term as

$$\|\vec{u}-\vec{v}\|_1^2 = (|-2|+|-3|+|4|+|2|)^2 = (2+3+4+2)^2 = 11^2 = 121.$$

The second term of Eq. A8 evaluates to

$$\|\vec{u}-\vec{v}\|_2^2 = (-2)^2 + (-3)^2 + 4^2 + 2^2 = 4+9+16+4 = 33.$$

The last term of Eq. A8 evaluates to

$$2\sum_{i=1}^{n-1}\left\{(|u_i-v_i|)\left(\sum_{j=i+1}^{n}|u_j-v_j|\right)\right\} = 2\left[|-2|(|-3|+|4|+|2|)+|-3|(|4|+|2|)+|4|(|2|)\right]$$
$$= 2\left[2(3+4+2)+3(4+2)+4(2)\right]$$
$$= 2\left[2(9)+3(6)+4(2)\right] = 2(18+18+8) = 88$$

Again, we see that the LHS evaluates to 121 and the RHS evaluates to $33+88=121$, supporting the validity of the relation presented in Eq. A8.

### Appendix B – Illustrating that the Joint Distance Measure is not a Distance Metric

The task is to check if given any three vectors $\vec{u}_0$ $\vec{u}_1$ and $\vec{u}_2$:

$$D_{JDM}(\vec{u}_0,\vec{u}_2) \leq D_{JDM}(\vec{u}_0,\vec{u}_1) + D_{JDM}(\vec{u}_1,\vec{u}_2). \tag{B1}$$

If the inequality in B1 holds, then the JDM is a distance metric. Here, we consider only the JDM 2 (combination of the Euclidean distance and the cosine similarity). Under this case, the inequality in B1 may be rewritten as:

$$\|\vec{u}_2-\vec{u}_0\|_2 |\cos\theta_{0,2}-2| \leq \|\vec{u}_1-\vec{u}_0\|_2 |\cos\theta_{0,1}-2| + \|\vec{u}_2-\vec{u}_1\|_2 |\cos\theta_{1,2}-2|, \tag{B2}$$

where $\theta_{0,1}$ is the angle between vectors $\vec{u}_0$ and $\vec{u}_1$, $\theta_{1,2}$ is the angle between vectors $\vec{u}_1$ and $\vec{u}_2$, $\theta_{0,2}$ is the angle between vectors $\vec{u}_0$ and $\vec{u}_2$. We note that in general, $\theta_{0,2}$ is not the sum of $\theta_{0,1}$ and $\theta_{1,2}$.

If we let $\vec{u}_1-\vec{u}_0 = \vec{v}$ and $\vec{u}_2-\vec{u}_1 = \vec{w}$ so that $\vec{u}_2-\vec{u}_0 = \vec{v}+\vec{w}$, then we may pose B2 as:

$$\|\vec{v}+\vec{w}\|_2 |\cos\theta_{0,2}-2| \leq \|\vec{v}\|_2 |\cos\theta_{0,1}-2| + \|\vec{w}\|_2 |\cos\theta_{1,2}-2|. \tag{B3}$$

To check if the inequality in B3 is valid, we start by taking the square of the LHS of B3:

$$\left( \|\vec{v} + \vec{w}\|_2 \left|\cos\theta_{0,2} - 2\right| \right)^2 = \|\vec{v} + \vec{w}\|_2^2 \left(\cos\theta_{0,2} - 2\right)^2. \tag{B4}$$

Using the fact that $\|\vec{v} + \vec{w}\|_2^2 = \|\vec{v}\|_2^2 + 2\vec{v}\cdot\vec{w} + \|\vec{w}\|_2^2$, we may expand B4 as:

$$\left( \|\vec{v} + \vec{w}\|_2 \left|\cos\theta_{0,2} - 2\right| \right)^2 = \left( \|\vec{v}\|_2^2 + 2\vec{v}\cdot\vec{w} + \|\vec{w}\|_2^2 \right) \left|\cos\theta_{0,2} - 2\right|^2. \tag{B5}$$

We expand the RHS of B5 to obtain:

$$\left( \|\vec{v} + \vec{w}\|_2 \left|\cos\theta_{0,2} - 2\right| \right)^2 = \left( \|\vec{v}\|_2 \left|\cos\theta_{0,2} - 2\right| \right)^2 + 2\vec{v}\cdot\vec{w}\left(\cos\theta_{0,2} - 2\right)^2 + \left( \|\vec{w}\|_2 \left|\cos\theta_{0,2} - 2\right| \right)^2. \tag{B6}$$

But by Cauchy-Schwarz relation [50, 51], $|\vec{v}\cdot\vec{w}| \leq \|\vec{v}\|_2^2 \|\vec{w}\|_2^2$.

Therefore, the equality relation in B6 can be replaced by the inequality:

$$\left( \|\vec{v} + \vec{w}\|_2 \left|\cos\theta_{0,2} - 2\right| \right)^2 \leq \left( \|\vec{v}\|_2 \left|\cos\theta_{0,2} - 2\right| \right)^2 + 2\|\vec{v}\|_2^2\|\vec{w}\|_2^2 \left(\cos\theta_{0,2} - 2\right)^2 + \left( \|\vec{w}\|_2 \left|\cos\theta_{0,2} - 2\right| \right)^2, \tag{B7}$$

The expression on the RHS of B7 is indeed a perfect square. Hence, B7 may be written as:

$$\left( \|\vec{v} + \vec{w}\|_2 \left|\cos\theta_{0,2} - 2\right| \right)^2 \leq \left( \|\vec{v}\|_2 \left|\cos\theta_{0,2} - 2\right| + \|\vec{w}\|_2 \left|\cos\theta_{0,2} - 2\right| \right)^2, \tag{B8}$$

By taking the square-root of both sides of B8, we obtain:

$$\|\vec{v} + \vec{w}\|_2 \left|\cos\theta_{0,2} - 2\right| \leq \|\vec{v}\|_2 \left|\cos\theta_{0,2} - 2\right| + \|\vec{w}\|_2 \left|\cos\theta_{0,2} - 2\right|, \tag{B9}$$

Now, comparing B9 to B3, we noticed that B3 will hold if $\cos\theta_{0,2} \leq \cos\theta_{0,1}$ and $\cos\theta_{0,2} \leq \cos\theta_{1,2}$. To check this, we consider two cases.

**Case 1 – $\theta_{0,2}$ is not the sum of $\theta_{0,1}$ and $\theta_{1,2}$**

In this case, $\theta_{0,2} \neq \theta_{0,1} + \theta_{1,2}$. Since there is no relation between these three angles, there is no guarantee that $\cos\theta_{0,2} \leq \cos\theta_{0,1}$ and $\cos\theta_{0,2} \leq \cos\theta_{1,2}$ will hold. Therefore, it is possible that for some values of $\theta_{0,1}$, $\theta_{1,2}$ and $\theta_{0,2}$ the inequality may hold and for some others, the relation may not hold.

**Case 2 – $\theta_{0,2}$ is the sum of $\theta_{0,1}$ and $\theta_{1,2}$**

In this case, $\theta_{0,2} = \theta_{0,1} + \theta_{1,2}$, and we may expand the cosine of this angle as:

$$\cos\left(\theta_{0,1} + \theta_{1,2}\right) = \cos\theta_{0,1} \cos\theta_{1,2} - \sin\theta_{0,1} \sin\theta_{1,2}. \tag{B10}$$

Thus,

$$\|\vec{v}\|_2 |\cos\theta_{0,2} - 2| = \|\vec{v}\|_2 |\cos\theta_{0,1} \cos\theta_{1,2} - \sin\theta_{0,1} \sin\theta_{1,2} - 2|. \tag{B11}$$

Since $-1 \leq \cos\theta \leq 1$ we know that $\theta_{0,1}$ and $\theta_{1,2}$ are constrained to be in the range $180° \geq \theta \geq 0°$. Thus, without loss of generality, we let $0° \leq \theta_{0,1} \leq \theta_{1,2} \leq 180°$, so that $\sin\theta_{0,1} \geq 0$, $\sin\theta_{1,2} \geq 0$ and thus, $\sin\theta_{0,1} \sin\theta_{1,2} \geq 0$. The sines of the angles are all nonnegative because the sine of any angle in the first or second quadrant is nonnegative. This simple fact implies that the term $\sin\theta_{0,1} \sin\theta_{1,2}$ in Eq. B11 is always nonnegative. Hence,

$$\|\vec{v}\|_2 |\cos\theta_{0,2} - 2| \leq \|\vec{v}\|_2 |\cos\theta_{0,1} \cos\theta_{1,2} - 2|. \tag{B12}$$

Similarly,

$$\|\vec{w}\|_2 |\cos\theta_{0,2} - 2| \leq \|\vec{w}\|_2 |\cos\theta_{0,1} \cos\theta_{1,2} - 2|. \tag{B13}$$

Applying B12 and B13 to B9, we may write:

$$\|\vec{v} + \vec{w}\|_2 |\cos\theta_{0,2} - 2| \leq \|\vec{v}\|_2 |\cos\theta_{0,1} \cos\theta_{1,2} - 2| + \|\vec{w}\|_2 |\cos\theta_{0,1} \cos\theta_{1,2} - 2|, \tag{B14}$$

We note that $|\cos\theta_{0,1} \cos\theta_{1,2} - 2| \leq |\cos\theta_{0,1} \cos\theta_{1,2}| + 2$.

Furthermore, we note the validity of the inequalities $|\cos\theta_{0,1} \cos\theta_{1,2}| + 2 \leq |\cos\theta_{0,1}| + 2$ and $|\cos\theta_{0,1} \cos\theta_{1,2}| + 2 \leq |\cos\theta_{1,2}| + 2$. Thus, we may write:

$$\|\vec{v} + \vec{w}\|_2 |\cos\theta_{0,2} - 2| \leq \|\vec{v}\|_2 (|\cos\theta_{0,1}| + 2) + \|\vec{w}\|_2 (|\cos\theta_{1,2}| + 2), \tag{B15}$$

A careful comparison of B15 to B3, shows that the RHS of B15 is larger than the RHS of B3. Thus, there is no guarantee that B3 will hold.